\documentclass[a4paper,10pt]{amsart}
\usepackage[utf8]{inputenc}
\title{Segments on the Right Branch of a Binary Tree}
\author{Naomi Lindenstrauss}
\email{Naomi.Lindenstrauss@mail.huji.ac.il}
\address{The Hebrew University of Jerusalem}
\date{December 18, 2020}
\begin{document}
\begin{abstract}
It is proved that the average number of segments on the right branch of a binary tree of size $n$ tends to $3$ as $n$ tends to $\infty$ . Also the fraction of trees with $k$ segments on the right branch from all trees of size $n$ tends to $\frac{k}{2^{k+1}}$ as $n$ tends to $\infty$ .
%
\end{abstract}
\maketitle

\section{Introduction}
A {\bf binary tree} is a structure defined recursively to be either a single {\em external node} or an {\em internal node} that is connected to two binary trees, a {\em left subtree} and a {\em 
right subtree}. The number of binary trees with $n+1$ external nodes is the Catalan number
\[ c_{n}  = \frac{1}{n+1} \left( \begin{array}{c} 2n \\ n \end{array} \right) \] 
 We call $n$ the size of the tree. (cf. \cite{fl} )

We prove that the average number of segments on the right branch of a binary tree of size $n$  (i.e. the sum of the numbers of segments on the right branch of all $n$ sized trees divided by $c_{n}$) tends to $3$ as $n$ tends to $\infty$. We also prove that for $n$ sized trees the relative part of trees with $k$ segments on the right branch  tends to 
\[\frac{k}{2^{k+1}}\]
as $n$ tends to $\infty$.

Both results were originally proved by generating functions methods.  Later a more combinatorial proof of the first result was found.

This paper is a byproduct of an attempt to investigate the average case complexity of the satisfiability probability of Boolean expressions, that can be represented as binary trees. The {\em construction step} in this paper can be used as an induction step to pass from $n$ sized Boolean expressions to $n+1$ sized Boolean expressions.

\section{The average number of segments on the right branch of a binary tree tends to $3$}
The proof is based on a construction of all $n+1$ binary trees from the $n$ sized binary trees. For each node on the right branch of an $n$ sized tree we create an $n+1$ sized tree which has at this node a subtree whose left subtree is what we had at this node and whose right subtree is an external node. For example

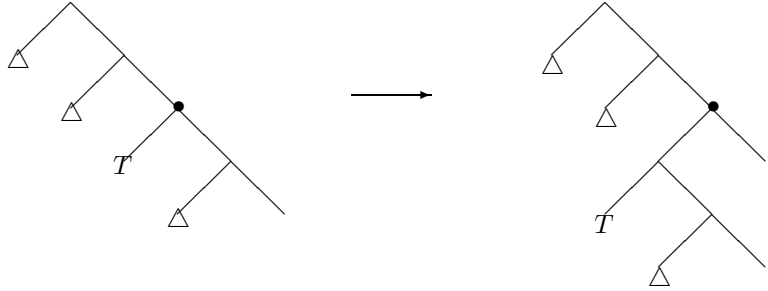
\begin{picture}(300,120)
\put(30,110){\line(1,-1){80}}
\put(30,110){\line(-1,-1){20}}
\put(50,90){\line(-1,-1){20}}
\put(70,70){\line(-1,-1){20}}
\put(90,50){\line(-1,-1){20}}
\put(68,68){$\bullet$}
\put(6,85){$\triangle$}
\put(26,65){$\triangle$}
\put(46,45){$T$}
\put(66,25){$\triangle$}
\put(135,75){\vector(1,0){30}}

\put(230,110){\line(1,-1){60}}
\put(230,110){\line(-1,-1){20}}
\put(250,90){\line(-1,-1){20}}
\put(270,70){\line(-1,-1){40}}
\put(250,50){\line(1,-1){40}}
\put(270,30){\line(-1,-1){20}}
\put(268,68){$\bullet$}
\put(206,83){$\triangle$}
\put(226,63){$\triangle$}
\put(226,23){$T$}
\put(246,3){$\triangle$}

\end{picture}

We have for example

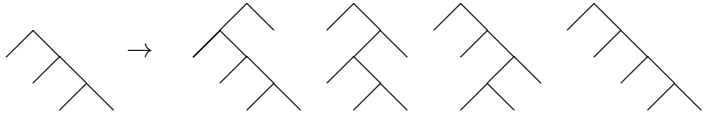
\begin{picture}(260,60)(0,-10)

\put(10,30){\line(1,-1){30}}
\put(10,30){\line(-1,-1){10}}
\put(20,20){\line(-1,-1){10}}
\put(30,10){\line(-1,-1){10}}
\put(45,20){$\rightarrow$}
\put(90,40){\line(-1,-1){20}}
\put(90,40){\line(1,-1){10}}
\put(80,30){\line(1,-1){30}}
\put(80,30){\line(-1,-1){10}}
\put(90,20){\line(-1,-1){10}}
\put(100,10){\line(-1,-1){10}}

\put(130,40){\line(-1,-1){10}}
\put(130,40){\line(1,-1){20}}
\put(140,30){\line(-1,-1){20}}
\put(130,20){\line(1,-1){20}}
\put(140,10){\line(-1,-1){10}}

\put(170,40){\line(1,-1){30}}
\put(170,40){\line(-1,-1){10}}
\put(180,30){\line(-1,-1){10}}
\put(190,20){\line(-1,-1){20}}
\put(180,10){\line(1,-1){10}}

\put(220,40){\line(1,-1){40}}
\put(220,40){\line(-1,-1){10}}
\put(230,30){\line(-1,-1){10}}
\put(240,20){\line(-1,-1){10}}
\put(250,10){\line(-1,-1){10}}

\end{picture}

Every $n+1$ sized tree is created in such a way. The subtree at the last but one node on the right branch has to be replaced by its left subtree. For example
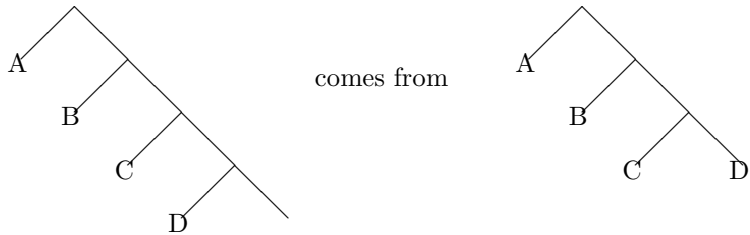
\begin {picture}(290,110)
\put(30,100){\line(1,-1){80}}
\put(30,100){\line(-1,-1){20}}
\put(50,80){\line(-1,-1){20}}
\put(70,60){\line(-1,-1){20}}
\put(90,40){\line(-1,-1){20}}
\put(120,70){comes from}
\put(220,100){\line(1,-1){60}}
\put(220,100){\line(-1,-1){20}}
\put(240,80){\line(-1,-1){20}}
\put(260,60){\line(-1,-1){20}}

\put(5,75){A}
\put(25,55){B}
\put(45,35){C}
\put(65,15){D}

\put(195,75){A}
\put(215,55){B}
\put(235,35){C}
\put(275,35){D}
\end{picture}

Using this construction step we can easily compute the average number of segments on the right branch for each $n$. Let us use the notation $i\times j$ to say that there are $i$ trees with $j$ segments on the right branch. 

For $n=1$ there is one tree that has one segment on the right branch, that is $1\times 1$,
and the average number of segments on the right branch is $\frac{1}{1} =1$.

For $n=2$ we have $1\times 1,\; 1\times 2$ so 
\[average = \frac{1*1+1*2}{1+1} =\frac{3}{2} =1.5 \]

For $n=3$ we have $2\times 1,\; 2\times2,\; 1\times 3$ so
\[average = \frac{2*1+2*2+1*3}{2+2+1}=\frac{9}{5}=1.8\]

For $n=4$ we have $5\times 1,\; 5\times 2,\; 3\times 3,\; 1\times 4$ so
\[average=\frac{5*1+5*2+3*3+1*4}{5+5+3+1}=\frac{28}{14}=2\]

For $n=5$ we have $14\times 1,\; 14\times 2,\; 9\times 3,\; 4\times 4,\;1\times 5$ so
\[average=\frac{90}{42}=2.14\]

For $n=6$ we have $42\times 1,\; 42\times 2,\; 28\times 3,\; 14\times 4,\; 5\times 5,\; 1\times 6$ so
\[average=\frac{297}{132}=2.25\]

For $n=7$ we have $132\times 1,\; 132\times 2,\; 90\times 3,\; 48\times 4,\;20\times5,\;6\times 6,\;1\times 7$ so
\[average=\frac{1001}{429}=2.33\]

For $n=8$ we have $429\times 1,\; 429\times 2,\; 297\times 3,\; 165\times 4,\; 75\times 5,\; 27\times 6,\; 7\times 7,\; 1\times 8$ so
\[average=\frac{3432}{1430}=2.4\]

For $n=9$ we have $1430\times 1,\;1430 \times 2,\; 1001\times 3,\;572\times 4,\; 275\times 5,\; 110\times 6.\; 35\times 7,\; 8\times 8,\; 1\times 9$ so
\[average=\frac{11934}{4862}=11934=2.45\]

for $n=10$ we have $ 4862\times 1.\; 4862\times 2,\; 3432\times 3,\; 2002\times 4,\; 1001\times 5,\; 429\times 6,\; 154\times 7,\; 44\times 8,\; 9\times 9,\; 1\times 10$
\[average=\frac{41990}{16796}=2.5\]

 I am  grateful to Michael Larsen for observing that in the expressions for the average, the numerator is the difference between consecutive Catalan numbers. It turns out that with the help of the construction step we can prove that the average number of segments on the right branch of a binary tree of size $n$ is
\[\frac{c_{n+1}-c_{n}}{c_{n}}\]

The proof is easy: From the construction step we see that the number of $n+1$ sized trees that an $n$ sized tree creates is 
\[1+the \; number\; of\; segments\; on\; its\; right\; branch.\]
So the sum of the numbers of the segments on the right branch of all $n$ sized trees is $c_{n+1}-c_{n}$ and the average number of the segments on the right branch of $n$ sized trees is
\[\frac{c_{n+1}-c_{n}}{c_{n}}\]
This tends to $3$ as $n$ tends to $\infty$ since it is equal to
\[\frac
{\frac{1}{n+2}\left( \begin{array}{c} 2n+2\\n+1 \end{array}\right) -\frac{1}{n+1} \left( \begin{array}{c} 2n\\n \end{array} \right)}
{\frac{1}{n+1}\left( \begin{array}{c} 2n\\n \end{array} \right) }  = \frac
{\frac{n+1}{n+2}(2n+1)(2n+1)}{(n+1)(n+1)} -1
\]

\section{The relative part of binary trees with $k$ segments on the right branch}
Let $B$ be the set of all binary trees. For a binary tree $T\epsilon B$ we denote by $|T|$ the number of its nodes (both internal and external). Let $N_{i}$ be the number of binary trees with $i$ nodes. We consider the generating function
\[N(z)=\sum_{i=1}^{\infty}N_{i}z^{i} =\sum_{T\epsilon B}z^{|T|} = \sum_{T=\,\bullet}z^{|T|} +\sum_{
\begin{picture}(40,40)
\put(0,30){$T=$}
 \put(20,30){$\bullet$}
 \put(22,30){\line(-1,-1){10}}
 \put(22,30){\line(1,-1){10}}
 \put(5,15){$T_{1}$}
 \put(33,15){$T_{2}$}
  \end{picture}} z^{|T|}  = \]
 \[=z+\sum_{T_{1},T_{2}\epsilon B}z^{1+|T_{1}|+|T_{2}|} 
 =z+z\sum_{T_{1}\epsilon B,T_{2}\epsilon B}z^{|T_{1}|}z^{|T_{2}|}=z+zN(z)^{2}
 \]

 Let $B_{k}$ be the set of all binary trees with $k$ segments on the right branch and  let $S_{i}^{k}$  be the number of binary trees of $i$ nodes that have $k$ segments on the right branch.
 Consider the generating function
 \[S^{k}(z) =\sum_{i=1}^{\infty}S_{i}^{k}z^{i} = \sum_{T\epsilon B_{k}} z^{|T|} =z^{k+1}N(z)^{k}
 \]

 Substitute 
 \[z=\frac{N}{1+N^{2}}\]
 in this formula obtaining
 \[S^{k}={\frac
 {N^{k+1}}
 {(1+N^{2})^{k+1} }N^{k}}=
 \frac{N^{2k+1}}{(1+N^{2})^{k+1}}\]
 (here we shortened $N(z)$  and $S^{k}(z)$ respectively to $N$ and $S^{k}$).
 
 We have for $N(z)$ an equation of the form $N(z)=z\:\Phi (N(z))$ (in our case $\Phi (x)=1+x^{2}$). Let $\tau$ be the smallest positive root of the equation
 \[\Phi (x)=x\:\Phi '(x)\]
 (in our case $\tau =1$). 
 
 From the methods outlined in \cite{lind,kap} it follows that  
 \[\lim  _{i\rightarrow \infty}\left.\frac{S_{i}^{k}}{N_{i}}=\frac{dS^k}{dN}
 \right|_{N=\tau}\]
 hence in our case
 \[\lim _{i \rightarrow \infty}\frac{S_{i}^{k}}{N_{i}}=\frac{d}{dN}\left.\left(\frac{N^{2k+1}}{(1+N^{2})^{k+1}}\right)\right|
 _{N=1}=\]
\[ =\left.\frac{(1+N^{2})^{k+1}(2k+1)N^{2k}
 -N^{2k+1}(k+1)(1+N^{2})^{k}2N}
 {(1+N^{2})^{2k+2}}\right|_{N=1}=\frac{k}{2^{k+1}}
\]
In other words the fraction of binary trees with  $k$ segments on the right branch from all binary trees of size $n$ tends to $\frac{k}{2^{k+1}}$ as $n$ tends to $\infty$.

As can be expected \[\sum_{k=1}^{\infty} \frac{k}{2^{k+1}} =1\;\;\;\; and\;\;\;\; \sum_{k=1}^{\infty}\frac{k^{2}}{2^{k+1}}=3\]
These results can be easily proved by differentiating
\[x+x^2+x^3+\cdots = \frac{x}{1-x}\]

\end{document}